\documentclass{article}
\usepackage{amsbsy,amssymb,amscd,amsfonts,latexsym,amstext,delarray,
amsmath, diagrams}  \headsep=15pt
\def\cP{ {\cal P} }

\def\cX{{\cal X} }

\def\bX{ \bar{X}}
\def\Ext{{ \mbox{Ext} }}
\def\rank{{ \mbox{rank} }}
\def\dim{{ \mbox{dim} }}
\def\Spec{{ \mbox{Spec} }}

\def\Hom{{ \mbox{Hom} }}

\def\om{{ \omega  }}
\def\O{{ {\cal O} }}

\def\ra{{ \rightarrow }}

\def\g{{ \gamma }}
\def\cR{ {\cal R}}
\def\d{{ \delta }}

\def\e{{ \epsilon }}
\def\F{ {\mathbb F} }
\def\hZ{ \hat{\Z}}

\def\hra{{ \hookrightarrow }}

\def\C{{ \mathbb{C} }}
\def\R{{ \mathbb{R}}}

\def\G{{ \Gamma }}
\def\Gal{{ \mbox{Gal} }}

\def\barF{ \bar{\mathbb{F}}}
\def\bQ{\bar{\Q}}

\def\cE{ {\cal E}}

\def\Z{{ \mathbb{Z}}}

\def\s{ \sigma}

\usepackage{amsmath}
\usepackage{amsfonts}
\usepackage{amscd}
\usepackage{amssymb}
\newarrow{Corresponds}<--->
\def\Q{\mathbb{Q}}

\def\invlim{\varprojlim}

\def\ord{\mbox{ord}}

\def\P{ {\bf P}}
\def\cY{{\cal Y}}
\def\bq{\begin{quote}}
\def\eq{\end{quote}}

\def\all{\forall}

\def\cE{ {\cal E}}
\def\D{\Delta}
\title{Classical Motives : Motivic $L$-functions}
\author{Minhyong Kim}
\begin{document}

\maketitle
The exposition here follows the lecture delivered at the summer school, and hence,
contains neither precision, breadth of comprehension, nor depth of insight. The goal rather
 is the curious one of providing a loose introduction to the excellent introductions that
already exist, together with scattered parenthetical commentary.
The inadequate nature of the exposition is certainly worst in the third
section. As a remedy, the article of Schneider \cite{schneider} is recommended as a good starting point for the
complete novice, and  that of Nekovar \cite{nekovar} might be consulted for more streamlined formalism.
For the Bloch-Kato conjectures, the paper of Fontaine and Perrin-Riou \cite{FP} contains  a very systematic treatment,
while Kato \cite{kato1} is certainly hard to surpass for inspiration.  Kings \cite{kings}, on the other hand,
gives a nice summary of results (up to 2003).

\section{Motivation}

Given a variety
 $X$ over $\Q$, it is hoped that
a suitable analytic function $$\zeta(X,s),$$
a  {\em $\zeta$-function}  of $X$,
 encodes important arithmetic invariants of $X$.
The terminology of course stems from the fundamental
function
$$\zeta(\Q,s)=\Sigma_0^{\infty}n^{-s}$$
 named by Riemann, which is interpreted in this
general context as the zeta function of $\Spec(\Q)$.
A general zeta function should generalize Riemann's function in
a manner similar to Dedekind's extension to number fields.
Recall that the latter can be defined by replacing the
sum over positive integers by a sum over ideals:
$$\zeta(F,s)=\Sigma_{I}N(I)^{-s}$$
where $I$ runs over the ideals of the ring of integers
$\O_F$ and $N(I)=|\O_F/I|$,
and that
$\zeta(F,s)$ has a simple pole at $s=1$ (corresponding to the trivial
motive factor of $\Spec(F)$, as it turns out) with
$$(s-1)\zeta(F,s)|_{s=1}=\frac{2^{r_1}(2\pi)^{r_2}h_F R_F}{w_F\sqrt{|D_F|}}$$
By the unique factorization of ideals, $\zeta(F,s)$ can also
be written as an Euler product
$$\prod_{\cP}(1-N(\cP)^{-s})^{-1}$$
as $\cP$ runs over the maximal ideals of $\O_F$, that is,
the closed points of $\Spec(\O_F)$.
Now, if a scheme $\cY$  is of finite type over $\Z$,
then for any closed point $y\in \cY$, its residue
field $k(y)$ is finite. Write $N(y):=|k(y)|$.
We can then form an Euler product \cite{serre1}
$$Z(\cY,s):=\prod_{y\in \cY_0}(1-N(y)^{-s})^{-1},$$
where $(\cdot)_0$ denotes the set of closed points for
any scheme $(\cdot)$. In the case when the map
$$\cY \ra \Spec(\Z)$$ factors through $\Spec(\F_p)$,
$Z(\cY,s)$  reduces to Weil's zeta function for a variety over
a finite field (with the substitution $p^{-s}\mapsto t$ if a formal variable
has intervened as in \cite{serre1}, section 1.6).

When we are starting with $X/\Q$,  a straightforward imitation of Dedekind's definition might involve
taking an integral model
$\cX$ of $X$, which is a proper flat scheme of finite-type
over $\Z$ with $X$ as generic fiber, and defining
$$\zeta(X,s)``:="Z(\cX,s)=\prod_{x\in \cX_0}(1-N(x)^{-s})^{-1}$$
The problem with this approach is that the function thus obtained
will depend on the model, and there is no general method for
choosing a canonical one. However, there will be some set $S$ of primes
such that there is a model $\cX_S$ over $\Spec(\Z[1/S])$
which is furthermore {\em smooth}. Even though such a $\Z[1/S]$-model need
be no more canonical, it does turn out that the incomplete zeta function
$$\zeta_S(X,s):=\prod_{x\in (\cX_S)_0}(1-N(x)^{-s})^{-1}$$
is independent of the model. (More on this point below.)
So there are good elementary generalizations of
incomplete zeta functions. We note in this connection that
$$Z(\cX,s)=\prod_pZ(\cX_p,s)$$
where
$$\cX_p=\cX\otimes \F_p$$
is the reduction of $\cX$ modulo $p$, so that that
$$\zeta_S(X,s)=\prod_{p\notin S} Z(\cX_p,s)$$
is the result of deleting a few Euler factors. Thus, the problem of
defining a canonical zeta function becomes one of
inserting canonical factors for the primes of bad reduction.
It is not impossible that there is a theory of integrals models
that isolates a class that is canonical enough to yield
a good definition of $\zeta(X,s)$. But the
current approach proceeds instead to
break up partial zeta function into natural factors
$$\zeta_S(X,s)=\prod L_S(M_i,s)^{\pm 1},$$
 according to the way $X$ is decomposed into
constituent  motives $\{ M_i\}$ in a suitable category.
(It is not much of an exaggeration to say that the decomposition of  zeta functions
is the main empirical phenomenon leading to the
hypothesis of a category of motives.)
The incomplete $L$-functions
 $L_S(M_i,s)$ of the $M_i$ should
then encode arithmetic invariants of the $M_i$,
which, in turn, refine the arithmetic invariants of
$X$. It is believed that good analytic properties
must be established to access the invariants efficiently,
including functional equations.
This, in turn, requires us to complete the $L$-functions
using cohomological machinery in general. The completed $L$-functions
then will lead to a completed zeta function.

A simple illustration is provided by
the elementary example of an elliptic curve
$E/\Q$  with affine  equation
$$y^2+a_1xy+a_3y=x^3+a_2x^2+a_4x+a_6$$

Let $\cE_S$ be a smooth and proper $\Z[1/S]$ model. Then
$$\zeta_S(E,s):=Z(\cE_S,s)$$
It is not very hard to check that
$$\zeta_S(E,s)=\zeta_S(\Q,s)\zeta_S(\Q,s-1)/L_S(H^1(E),s)$$
(\cite{silverman1}, V.2.4) illustrating the kind of decomposition alluded to
above. Here
$$\zeta_S(\Q,s)=\Sigma_{\{(n,p)=1, \all p\in S\}} n^{-s}$$
is a standard incomplete zeta function and
$$L_S(H^1(E),s)=\prod_{p\notin S} L_{p}(H^1(E),s)$$
is the incomplete $L$-function of $E$ with factors defined by
$$L_{p}(E,s)=\frac{1}{1-a_{p}p^{-s}+p^{1-2s}}$$
Here
$a_p=p+1-N_p$
and $N_p$ is the number of points on $E$ mod $p$.
$L_S(E,s)$ turns out to be the partial $L$-function
corresponding to the motivic factor $H^1(E)$ of $E$.

We can put in Euler factors for $p\in S$. It is
obvious how to do it for  $\zeta_S(\Q, s)$ and
$\zeta_S(\Q, s-1)$ giving us the Riemann
zeta function
$\zeta(\Q, s)$ and its shift $\zeta(\Q, s-1)$ respectively.
For the incomplete $L_S(H^1(E),s) $, we put in the factors
 according to a recipe determined by the
reduction of $E$ at $p$:
$$L_p(H^1(E), s)=\left\{ \begin{array}{cl}1/(1-p^{-s}) & \mbox{split multiplicative};\\
1/(1+p^{-s})& \mbox{non-split multiplicative};\\
1 & \mbox{additive}.\\
\end{array} \right.
$$
(\cite{silverman2}, II.10)
and define
$$L(H^1(E),s):=\prod_pL_p(H^1(E),s)$$
Here we have used  the breakdown of the incomplete zeta function
into three factors as an aid in defining the full zeta
function of $E$.
However, this case is somewhat misleading in that there {\em is} a
canonical model that could have been used instead, namely,
the Weierstrass minimal  model
$$\cE$$that appears in basic textbooks.
In fact, one can check that
$$\zeta(E,s)=\zeta(\Q,s)\zeta(\Q,s-1)/L(H^1(E),s)=Z(\cE,s)$$
as follows from the trace formula (\cite{silverman1}, V.2) for the Frobenius map on
elliptic curves for $p\notin S$, and a much easier
counting argument for $p\in S$. So this would seem to
be an instance where the naive extension of Dedekind's method
works out.
Nevertheless, we explain how the bad factors can be obtained without
reference to the model, starting at this point to use the language of
\'etale cohomology \cite{milne}.
In the sequel, we  fix an algebraic closure $\bQ$ of $\Q$,
closures $\bQ_p$ of $\Q_p$, and embeddings
$\bQ\hra \bQ_p$. Therefore, we have embeddings of
Galois groups
$$G_p:=\Gal(\bQ_p/\Q_p)\hra G:=\Gal(\bQ/\Q)$$
The residue field of $\bQ_p$ is an algebraic closure $\barF_p$ of
$\F_p$, and we have an exact sequence
$$0\ra I_p \ra G_p \ra \Gal(\barF_p/\F_p)\ra 0$$
defining the inertia subgroup $I_p$.
Denote by $Fr_p$  the generator of
$\Gal(\barF_p/\F_p)$ that takes $x$ to $x^{1/p}$. Finally,
$\bar{\cE_p}$ denotes the basechange of $\cE_p$ to
$\barF_p$ and $\bar{E}$ the basechange of $E$ to $\bQ$.
 We need the \'etale cohomology
 $$H^1(\bar{E},\Q_l)$$
 for primes $l$, and
 $$H^1(\bar{\cE_p},\Q_l)$$
 for $l\neq p$.
By the Lefschetz trace formula (\cite{milne}, VI.12.3),
$$Z(\cE_p,s)=\frac{\det ([I-p^{-s}Fr_p]|H^1(\bar{\cE_p},\Q_l))}{\det ([I-p^{-s}Fr_p]|H^0(\bar{\cE_p},\Q_l))\det ([I-p^{-s}Fr_p]|H^2(\bar{\cE_p},\Q_l))}$$
But for each $i=0,1,2$,
$$H^i(\bar{\cE_p},\Q_l)\simeq H^i(\bar{E},\Q_l)^{I_p}$$
the superscript referring to the subspace of elements fixed by the inertia action.
(For $H^0$ and $H^2$, this is an easy exercise.
The $H^1$ case  is slightly harder. See \cite{milne}, proof of theorem V.3.5. Although the discussion
there is given for
smooth surfaces fibered over `geometric' curves, it is rather straightforward to adapt it to
the present situation.)
For $p\notin S$, any pair $X\hra \cX $ as above satisfies
$$H^i(\bar{\cX_p},\Q_l)\simeq H^i(\bar{X},\Q_l)$$
where the $I_p$-action must be trivial,
and provides the reason that
the incomplete zeta function is independent of the model (\cite{milne}, VI.4.1).
In any case, it ends up  that the bad factor could have been written
$$Z(\cE_p,s)=\frac{\det ([I-p^{-s}Fr_p]|H^1(\bar{E},\Q_l)^{I_p})}{\det ([I-p^{-s}Fr_p]|H^0(\bar{E},\Q_l)^{I_p})\det ([I-p^{-s}Fr_p]|H^2(\bar{E},\Q_l)^{I_p})}$$
in a way that refers only to $E$.
It is this formula that generalizes to arbitrary motives.

Since we have thus far been entirely cavalier about convergence, we note in passing that
the estimate $|a_p|\leq 2 \sqrt{p}$ (\cite{silverman1}, V.II)
implies that the Euler product converges for
$Re(s)> 3/2$.

To control fine analytic properties, one establishes a relation to
automorphic $L$-functions.
For elliptic curves such a relation can be made explicit by computing the conductor
$$N_E:=\prod_{p\in S}p^{f_p}$$
Here
$$f_p=ord_p(\D_E)+1-m_E$$
where $\D_E$ is the discriminant of $E$ and $m_E$ is the number of geometric components (that is,
components over $\bar{\F}_p$) of the special fiber of the Neron model of $E$.
Even though this formula for $f_p$ again refers to the model, it can be defined purely
in terms of the Galois action on $H^1(\bar{E},\Q_l)$ (\cite{silverman2}, IV.10).

The well-known and deep
 fact, established through the work of Wiles, Taylor-Wiles, and Breuil-Conrad-Diamond-Taylor (\cite{wiles},
 \cite{TW}, \cite{BCDT}), is that
$L$ has an analytic continuation to the complex plane.
More precisely,
$$L(E,s)=L(f_E,s)=\frac{1}{(2\pi)^s\G(s)}\int_0^{\infty} f_E(iy)y^{s-1}dy$$
$$=\frac{1}{(2\pi)^s\G(s)}[\int_{1/\sqrt{N_E}}^{\infty}f_E(iy)y^{s-1}dy\pm \int_{1/\sqrt{N_E}}^{\infty}f_E(iy)y^{1-s}dy]$$
for a normalized weight 2 new cusp form $f_E$ of level $N_E$ which is an
eigenvector for the Hecke operators, determined by a $q$ expansion
$$f_E=q+a_2q^2+\cdots$$
where the  $a_p$ have to be the same as those for $E$ and the general coefficient
is determined by those with prime index.
Because the weight 2 cusp forms of level $N_E$ form a finite dimensional
space, it is easy to see that  $f_E$ is completely determined by a finite computation
of the $a_p$'s, and that this integral formula can then be used to compute
$L$-values.

The celebrated conjecture of Birch and Swinnerton-dyer (BSD) \cite{BSD} says
$$\ord_{s=1} L (E,s) =\rank E(\Q)$$
The equality is known if $\ord_{s=1}L(E,s) \leq 1$ by the work of Gross-Zagier and Kolyvagin
(\cite{GZ}, \cite{kolyvagin}).
Recall that
$$E(\Q)=\Z^r\times E(\Q)_{tor}$$
where the finite abelian group $E(\Q)_{tor}$ is easily computed using the Nagell-Lutz theorem
(\cite{silverman1}, VIII.7).
This conjecture promises to give an analytic approach to
understanding the elusive rank $r$.
However, it must be admitted that even though the $L$-function is computable, the utility of this equality
for actually computing the rank of an elliptic curve is somewhat ambiguous.
This is because the order of zero of an analytic function
might not be possible to determine using a finite computation.
We will discuss below how the vanishing itself can be computably
determined using the refined version of this conjecture.
On the other hand, an extremely useful viewpoint on the order of vanishing arises from the
functional equation. That is to say, one inserts a {\em gamma factor}, determined by the
Hodge theory of $E$ and viewed as the
contribution of the prime at infinity. With another correction factor contributed by the
conductor,
we arrive at a further completion:
$$\Lambda(E,s):=(2\pi)^s\G(s) N_E^{s/2}L(E,s)$$
which then satisfies a functional equation
$$\Lambda(E,2-s)=w_E\Lambda (E,s)$$
where $w_E=\pm 1$ depends on the curve $E$.
In fact, $w_E$ can be expressed as a product of
local terms
$$w_E=\prod_pw_{E,p}$$
each of  which can be computed in a straightforward fashion (\cite{BSD2}, section 6).

A significant corollary is that the {\em parity} of the
order is determined by the sign of $w_E$, usually referred to
as the  {\em sign of the functional equation}.
For example,  if $w_E=-1$, then clearly
$$L(E,1)=0$$
Suppose you can check
$L'(E,1)\neq 0$ using the equality with
$L'(f,1)$ (non-vanishing {\em can} be verified!). Then
we conclude that $E(\Q)$ has rank one.
Thus, one can produce many examples where the refined analysis of the L-function, including
the functional equation and computation, gives us the complete structure of $E(\Q)$.
(For deeper developments in this direction, see \cite{MR}.)

The BSD conjecture  continues as follows. If $r$ is the order of vanishing, then
$$(s-1)^{-r}L(E,s)|_{s=1}=|Sha(E)|R_E\Omega \prod_{p} c_p/|E(\Q)_{tor}|^2.$$
That is to say, the $L$-function purportedly encodes refined Diophantine
invariants of $E$, which we proceed to describe briefly (\cite{silverman1}, C.16).

The set $E(\Q)_{tor}$, which occurred already above, is the (finite) torsion subgroup of $E(\Q)$.
The $c_p$ refer to  the  Tamagawa numbers at primes $p$, consisting of the index
$$c_p=(E(\Q_p):E^0(\Q_p))$$
where $E^0(\Q_p)\subset E(\Q_p)$
is the subset of points that reduce to the connected component of the identity
in the Neron model of $E$ (\cite{silverman1}, VII.2). In particular, $c_p=1$ for primes of good reduction.
The difficult rational term is the order of
$Sha(E)$, the Tate-Shafarevich group of $E$, conjectured to be finite.
It is defined as the kernel
$$0\ra Sha(E) \ra H^1(G,E(\bQ))\ra \prod_p H^1(G_p,E(\bQ))$$
of the localization map on classifying spaces of torsors for $E$ in the \'etale topology of
$\Spec(\Q)$
as $p$ runs over all primes of $\Q$.

Then there are the transcendental terms:
$\Omega$, the real period (or twice that), defined as an integral
$$\Omega=\int_{E(\R)}|\om|$$
where $\om=dx/(2y+a_1x+a_3)$. The period can be easily
computed, but the inaccessible part is the regulator
$R_E$.  This is the covolume of
the lattice $E(\Q)/E(\Q)_{tor}$ inside the
inner product space
$$([E(\Q)/E(\Q)_{tor}]\otimes \R, <\cdot, \cdot>)$$
where  $< \cdot, \cdot >$ is defined by the
Neron-Tate canonical height. Thus, if
$\{P_1, P_2,\ldots, P_r\}$ is a basis for
for $E(\Q)/E(\Q)(tor)$, then
$$R_E:=|\det(<P_i,P_j>)|$$
Obviously, computation of $R_E$ would require knowledge of the
Mordell-Weil group. On the other hand, since the  formula gives a computable
bound for the denominator of $L(E,1)/\Omega$ when $R_E=1$, assuming its validity
allows us to verify the vanishing of $L(H^1(E),1)$ after a finite computation.
We refer the reader to \cite{RS} for an accessible report
 on the BSD conjecture, covering work up to 2002.

The known relations between $L$-functions and arithmetic are expected to
generalize vastly.
As indicated above, $L$-functions are defined using Galois actions on \'etale cohomology and
completed using Hodge theory.
\medskip

Before we summarize the relevant definitions, we recall the
big picture represented by the following conjectures:

(1) Hasse-Weil: the completed $L$ -function has a meromorphic continuation to the
complex plane and satifies a functional equation. This conjecture is supposed to
be addressed
by Langlands' program, which says `Motivic L-functions are automorphic L-functions'
( \cite{langlands}).

(2) Conjectures about values:

(a) Deligne's conjecture gives a general definition  of periods (in non-vanishing case)
using comparison of rational De Rham and topological cohomologies (\cite{deligne1});

(b) The Beilinson conjectures continue the discussion of orders of vanishing  and the
regulator using the rank and covolume of motivic cohomology (\cite{beilinson1}, \cite{beilinson2}, \cite{beilinson3},
\cite{beilinson4})

(c) The Bloch-Kato conjectures generalize the discussion of the rational part using
Tamagawa numbers (or determinants) for Galois representations via $p$-adic Hodge theory
(\cite{BK}, \cite{kato1}, \cite{FP}).
\medskip

\section{Definitions}
Let $X/\Q$ be a smooth projective variety as before.
Associated to $X$ is a well-known collection of cohomology
groups, the {\em realizations} of the motive(s) of $X$.
\bq
$H^n_{l}(X)=H^n_{et}(\bX, \Q_l)$ for each prime $l$, the $\Q_l$-coefficient \'etale cohomology
of degree $n$. This carries a natural action
of $G=\Gal(\bQ/\Q)$.

$H^n_{DR}(X):=H^n(X , \Omega_X^.)$, the algebraic De Rham
cohomology equipped with a Hodge filtration given by
$$F^iH^n_{DR}(X)=H^n(X, \Omega^{\geq i})\hra H^n_{DR}(X)$$
for each $i$.

$H^n_B(X):=H^n(X(\C), \Q)$, the $\Q$-coefficient singular cohomology of the
complex manifold $X(\C)$ equipped with a continuous
action $F_{\infty}$ of complex conjugation.\eq
The completed $L$-function of $H^n(X)$ uses all these structures.
These cohomology groups  are bound together by an intricate system of
canonical comparison isomorphisms. For example,
$$H^n_B(X)\otimes \Q_l\simeq H^n_{l}(X)$$
preserving the action of $F_{\infty}$, the complex conjugation. And then,
$$H^n_B(X)\otimes \C \simeq H^n_{DR}(X)\otimes \C.$$
This isomorphism endows the pair $(H^n_B(X), H^n_{DR}(X)\otimes \R)$ with
a {\em rational Hodge structure} of weight $n$ `defined over $\R$.'
That is, we have a direct sum decomposition
$$H^n_B(X)\otimes \C\simeq \oplus H^{p,q}(X)$$
where
$$H^{p,q}:=F^p\cap \bar{F}^{q}$$
and
$$F_{\infty}(H^{p,q})=H^{q,p}$$
If we denote by $c$ the
complex conjugation on $\C$
then
$$(H^n_B(X)\otimes \C)^{F_\infty\otimes c}=H^n_{DR}\otimes \R$$
At non-archimedean places, there is  an important analogue.
For any prime $p$,  we have
$$D_{DR}(H^n_{p}(X)):=(H^n_{p}(X)\otimes B_{DR})^{G_p}\simeq H^n_{DR}(X)\otimes \Q_p$$
where
$B_{DR}$ is Fontaine's ring of $p$-adic periods \cite{fontaine}.

These structures taken together motivate the following observation.
Regardless of its precise definition, a motive $M$
should have associated to it a collection of objects
as above that we call {\em a pure system of realizations}
that make up a category $\cR$.
This is a collection
$$R(M)=\{\{M_{l}\}, M_{DR}, M_B\}$$
where each $M_{l}$ is a representation of $G$ on a (finite-dimensional)
$\Q_l$-vector space, $M_{DR}$ is a filtered $\Q$-vector space,
and $M_B$ is a $\Q$-vector space with an involution $F_{\infty}$.
These vector spaces should all have the same
dimension and  be equipped with a system of comparison
isomorphisms as above.
The data must be subject to further subtle constraints
having to do with local Galois representations.

That is to say
recall the exact sequence:
$$0\ra I_p \ra G_p \stackrel{v}{\ra}  \Gal (\bar{\F}_p/ \F_p)\ra 0$$
For $l\neq p$,
$I_p$ has a tame $l$-quotient
$$t_l:I_p\ra I_{p,l}$$with the structure
$$I_{p,l}\simeq \hZ_l(1)\simeq  \invlim \mu_{l^n}$$
as a module for $\Gal (\bar{\F}_p/ \F_p)$.
Define $$W_p:=v^{-1}(\Z)\subset G_p,$$
 the {\em Weil group} at $p$.
It is convenient to analyze the data of $M_l$ using an associated
{\em Weil-Deligne} (W-D) representation \cite{tate1}
$$WD_p (M_l)$$ for each $p$, consisting of
a  representation $r$ of
$W_p$ such that $r|I_p$ has finite image and a nilpotent operator $N_p$ acting
on the representation.

These satisfy a compatibility condition
$$r(\phi_p)N_p r(\phi_p^{-1})=p^{-1}N_p$$
for any lift $\phi_p\in W_p$ of $Fr_p$.

The construction of $WD_p(M_l)$ for $p\neq l$
 uses the fact that
the action of $G_p$ when restricted to
 some finite index subgroup $G'_p$
is semi-stable, i.e., its inertia subgroup $I'_p$ acts
unipotently. Hence, the action of $I'_p$
can be expressed as
$$\s \mapsto \exp (t_l(\s) N_p)$$
for a nilpotent $N_p$. Then the representation $r$ on $I_p$
is given by
$$r(\phi_p^n\s)=\phi_p^n\s \exp(-t_l(\s)N_p)$$
for some choice of $\phi_p$. In fact, since the data $(\phi_p, N_p)$
determine the Weil-Deligne representation, it is usual to
identify the representation with the such a pair.

For $p=l$, we use the fact that any De Rham representation
is potentially semistable \cite{colmez}, and hence, gives us
a filtered $(\phi_p, N_p)$ module via
$$M_p \mapsto (M_p\otimes B_{st})^{G'_p}$$
If $G'_p=G_p$ (that is, if the representation itself is
semi-stable), then this gives us a Weil-Deligne representation
in an obvious way by defining
$$r(g)=\phi_p^n$$
if $g\mapsto Fr_p^n\in \Gal (\bar{\F}_p/ \F_p)$.
In \cite{FP}, I.1, it is explained in detail how
one extracts a Weil-Deligne representation from the data in the general case.

The viewpoint of the Weil-Deligne representation allows us to parametrize the
information of the Galois representations in a form that
does not use the topology of $\Q_l$. It provides, thereby, a suitable framework
for comparing the representations as $l$ varies, and
makes natural the connection to
complex automorphic forms \cite{taylor}.
Furthermore, one
creates thereby a rather precise analogy with the theory of limit mixed Hodge structures
\cite{illusie}.

Now define the
{\em Frobenius semi-simplification} $WD_p(M_l)^{ss}$
of $WD_p(M_l)$ by  replacing $\phi_p$ with its semi-simple
part.
With the terminology thus introduced,
here are the constraints we impose
on our pure system of realizations:
\medskip

`Good reduction almost everywhere':
We assume  that there exists a finite set
$S$ of primes such that
$WD_p(M_l)$ is unramified for all $p\notin S$,
i.e.,
$N_p=0$ and $I_p$ acts trivially.

\medskip

`Algebraicity and independence of $l$':
There exists a Frobenius semi-simple W-D representation
$WD_p(M)$ over $\bQ$ such that
$$WD_p(M)\otimes \bQ_l\simeq WD^{ss}_p(M_l)\otimes \bQ_l$$
\medskip

Subject to these conditions,
the collection $\{M_l\}$ is then referred to
as a
{\em strongly compatible} system of $l$-adic representations.
\medskip

`Weil conjecture':
There should exist
an integer $n$, called the {\em weight} of $M$,
such that the eigenvavlues of $Fr_p$ acting on
$WD_p(M)$ for $p\notin S$ have all Archimedean
absolute values equal to $p^{n/2}$. Furthermore, the
Hodge structure $M_B$ should be pure of weight $n$.
\medskip

`Purity of the monodromy filtration':
If we denote by $Mn_.$ the unique increasing filtration on $WD_p(M)$
such that $Mn_{-k}=0$, $Mn_{k}=WD_p(M)$ for sufficiently large
$k$, and $$N(Mn_k)\subset Mn_{k-2},$$ then the associated graded
piece
$$Gr^{Mn}_k(WD_p(M))$$
has all Frobenius eigenvalues of archimedean absolute value
$p^{(n+k)/2}$.
\medskip

It should be remarked that in general, we need to allow coefficients in
$E_{\lambda}$ for the representations
where $E$ is a number field and $E_{\lambda}$
are completions. Such coefficient systems arise naturally when
 considering direct summands of $\Q_l$ representations
or {\em motives with extra endomorphisms},
e.g., abelian varieties with CM. We will omit this generality
in this summary.
Another interesting view is  that the bi-grading
$$M_{B}\otimes  \C\simeq \oplus M^{p,q},$$
which is compatible with the complex conjugation of coefficients,
corresponds to a representation of
the group
$$Res^{\C}_{\R}({\bf G}_m).$$
Together with the action of
$$F_{\infty}\circ C$$ it can be viewed as a representation of
the real Weil group $W_{\R}$ (cite{tate1}) with points given by
$$W_{\R}(\R)=\C^*\cup \C^*j$$
where
$j^2=-1$ and $jzj^{-1}=\bar{z}$.
Here, $C$ is the Weil operator defined by
$$C|M^{pq}=i^{q-p}$$

It is conjectured that the realizations
$$H^n(X)=(\{H^n_l(X)\}, H^n_B(X), H^n_{DR}(X))$$ coming from a smooth projective variety $X$
satisfy the algebraicity, independence of $l$, and purity
conditions even  for $p\in S$.

The category of pure motives should be comprised of
objects in $\cR$ of {\em geometric origin}, a notion
without an entirely precise interpretation \cite{deligne2}. For example, we
need to admit  at every stage duals (homology) and tensor
products of all objects considered.
Objects that are not generated in an obvious way from
those of the form
$$H^n(X)$$
arise via images (or kernels) under pull-backs
and push-forwards in cohomology induced
by  maps of varieties,
as well as  $\Q$-linear combinations of geometric maps.
We should also be able to compose pull-backs with push-forwards.
Such compositions give rise to
the idea of using {\em correspondences} modulo
homological equivalence as morphisms \cite{deligne1}.
Once morphisms are constructed in this manner, we
naturally obtain new objects
using the decomposition of
$$End(H^n(X)),$$
which is a semi-simple $\Q$-algebra subject to one of the
standard conjectures saying that numerical equivalence and
homological equivalence coincide \cite{jannsen2}.

One can consider also a category of mixed systems of realizations
by requiring a weight filtration
$$\cdots \subset  W_{n-1}M
\subset W_nM \subset W_{n+1}M \subset $$
compatible with all the comparisons
and such that
each graded quotient $$Gr^n_W(M)$$
is a pure system of realizations of weight $n$.
Mixed motives should be those of geometric origin
such as the cohomology of varieties that are not
necessarily smooth or proper.
But then, we need to include
objects like (finite-dimensional quotients of)
$$\Q[\pi_1]$$
or the (co)-homology of (co-)simplicial varieties \cite{jannsen}.

Given a pure system $M$ of realizations
 we can define its
$L$-function $L(M,s)$ as an Euler product
$$L(M,s)=\prod_pL_p(M,s)$$
with
$$L_p(M,s)=\frac{1}{\det [(1-p^{-s}Fr_p)|(WD_p(M))^{I_p=1, N_p=0}]}$$
If $M$ is of weight $n$, then the product converges (and hence is non-zero)
for $Re(s)>n/2+1.$ For some conceptual motivation for this definition based on duality in
the function field case, see \cite{deninger2}. (The point is that
the inertia fixed part is the stalk of the intermediate extension
of the \'etale sheaf corresponding to $M_l$.)

There is also a factor at $\infty$ depending upon the
representation $M_B\otimes \C$
of $W_{\R}$.
Define
$$\G_{\R}:=\pi^{-s/2}\G(s/2)$$
$$\G_{\C}:=2(2\pi)^{-s}\G(s)$$
$$h^{pq}:=\dim M^{pq}$$
$$h^{p,\pm}:=\dim M^{pp, \pm 1}$$
where the signs in the superscript refer to the $\pm 1$
eigenspaces of the $F_{\infty}$-action.
Then
$$L_{\infty}(M,s)$$
is defined as
$$
\prod_{p<q}\G_{\C}(s-p)^{h^{pq}} $$
for odd $n$, and$$
\prod_{p<q}\G_{\C}(s-p)^{h^{pq}}\G_{\R}(s-n/2)^{h^{n/2+}}
\G_{\R}(s-n/2+1)^{h^{n/2-}}$$
for $n$ even (\cite{serre2}).

It is conjectured that
$$\Lambda (M,s)=L_{\infty}(M,s)L(M,s)$$ has a meromorphic continuation to $\C$
and satisfies a functional equation
$$\Lambda(M,s)=\e(M,s)\Lambda(M^*,1-s)$$
where the epsilon factor has the form
$\e(M,s)=ba^s$. Note that the contribution of the conductor has also been incorporated into
this factor. (For a precise discussion of the factor in the case $M=H^n(X)$, see \cite{serre2}.)
As alluded to above, the general expectation is that this conjecture will be dealt with by the Langlands' program.

\section{Conjectures on zeros, poles, and values}

A list of overall references to this section should include the papers mentioned in the introduction
as well as the original articles \cite{bloch1}, \cite{deligne1}, \cite{beilinson1}, \cite{beilinson2},
\cite{beilinson3}, \cite{beilinson4}, and
\cite{BK}. Having mentioned thus the sources, we will then proceed to be somewhat sloppy with specific
citation.

Here is   some convenient notation:
\bq

$\Q$: trivial system of realizations.

$\Q(1):=H^2(\P^1)^*$

$\Q(i)=\Q(1)^{\otimes i}$
for $i\geq 0$

$\Q(i)=\Hom( \Q (-i),\Q)$
for $i<0$.
\eq
For a system $M$ of realizations we define its {\em Tate twists} by
tensor products with $\Q(1)$:
$$M(i):=M\otimes \Q(i)$$
Then for any smooth projective variety of dim $d$,
we have (\cite{milne}, VI.6)
$$H^{2d}(X)\simeq \Q(-d)$$
and a perfect pairing
$$H^i(X) \times H^{2d-i}(X) \ra H^{2d}(X)$$
Repeated cup product with the cohomology class of a hyperplane
gives us the hard Lefschetz theorem \cite{deligne3}
$$H^i(X) \simeq H^{2d-i}(X)(d-i).$$
The effect of the twisting on realizations is such that
$M(n)_{l}$ is the tensor product of $M_l$ with
the $n$-th power of the $\Q_l$ cyclotomic character giving the action
of $G$ on the $l$-power roots of unity, and
$$F^i(M(n)_{DR})=F^{n+i}M_{DR}$$
with a corresponding shift in Hodge numbers
$h^{pq}$.
Furthermore,
$$F_{\infty}|M(n)_B=(F_{\infty}|M_B)\otimes (-1)^n$$
and for the $L$-functions,
$$L(M(n),s)=L(M,s+n).$$
These facts all follow in an elementary way from the structure of
$H^2(\P^1)$.

To state the main conjectures without spending time on categorical
preliminaries,
we will focus on the case where $M$ is
$$H^n(X)=(\{H^n_l(X)\}, H^n_{DR}(X), H^n_B(X))$$ for a smooth projective variety $X$
of dimension $d$. Assume that $H^n(X)$ is a pure system of realizations so  that
the analytic continuation and functional
equation hold true.
Since we have
$$H^n(X)^*\simeq H^{2d-n}(X)(d)\simeq H^n(X)(n)$$
the functional equation relates
$$L(H^n(X),s)$$
and
$$L(H^n(X)(n), 1-s)=L(H^n(X),n+1-s)$$
with center of reflection
$$(n+1)/2.$$
 We will therefore confine interest mostly to
$$m\geq (n+1)/2$$ or, equivalently,
$$n+1-m\leq (n+1)/2,$$
and refer to $\mbox{Re}(s)=(n+1)/2$ as the
critical line. (The reader should consult Nekovar's elegant article
\cite{nekovar} for a careful discussion of how to relate the
points addressed here to those on the right of the critical line.)
In fact, it is conceptually
convenient to parametrize by the letter $m$ the integers
$\geq (n+1)/2$, but to focus  then on the values at
the points
$n+1-m$.
In the discussion of orders, however, we will make explicit the
case of $m=n/2$, and hence,
the possible pole at $n+1-m=n/2+1$ (just to the right of
the critical line), whose importance is evident from the classical example
of Dedekind's zeta functions.  As we will explain below, the general geometric importance
of this pole is related to Tate's conjecture on the cohomology classes of algebraic cycles
\cite{tate2}.
We remark also that the study of
$L(H^n(X), s)$ near $m$ corresponds to the study
of $L(H^n(X)(m),s)$ near $s=0$. Thus, when the conjectures are
formulated in terms of values at zero in the literature, one
encounters the assumption that the weight
$n-2m$ is negative.

We start then with the conjecture on orders.
The simple case arising from an elliptic curve was reviewed already where
 $n=1$ and $m=n+1-m$ is the reflection point $(n+1)/2=1$. The conjecture of Birch and
Swinnerton-dyer says
$$\ord_{s=1}L(H^1(E),s)=\rank E(\Q)$$
Now, an element  $$x\in E(\Q)$$ gives rise to an
extension in the
category $\cR$ of realizations
$$\d(x)\in Ext^1_{\cR} (\Q, H^1(E)(1))$$
via Kummer theory.
It is conjectured that when  $\cR$ is replaced by a suitable
category of motives, this is the only way to construct
such extensions. This notion conveys the basic flavor of  conjectures
on orders in the case of odd weight.

For an example in the even weight case, let $X=\Spec(F)$
for a number field $F$, assumed for simplicity to be Galois
over $\Q$. Then
$$\zeta(F,s)=L(H^0(\Spec(F)),s)$$
which breaks up as into a product of Artin $L$-functions
$$\zeta(F,s)\prod_{\rho}L(\rho,s)$$
as $\rho: \Gal(F/\Q)\ra Aut(V)$ runs over
finite-dimensional representations on $\Q$-vector spaces. This is the most basic example
 of a motivic decomposition.
 In fact, any finite-dimensional $\Q$-representation
defines an Artin motive $M_{\rho}$ and it is a theorem that
$$\ord_{s=1}L(\rho, s)=-\dim \Hom_{AM}(\Q, M_{\rho}),$$
where the Hom occurs inside the category of Artin motives \cite{deligne1}.
Note that $s=1$ in this case is the integer point $n/2+1$ just to the right of the
critical line mentioned above.

The general conjecture about orders  says
$$\ord_{s=n+1-m}L(H^n(X),s)$$
$$=\dim Ext^1_{Mot_{\Z}}(\Q, H^n(X)(m))-\dim \Hom_{ Mot_{\Z}} (\Q, H^n(X)(m))$$
The $\Hom$ and $\Ext$ should occur inside a conjectural category
of mixed motives over $\Z$ with $\Q$-coefficients.
For weight reasons, the $\Hom$ term should vanishes unless
$n=2m$ in which case the $\Ext $ term should vanishes. That is,
in the pure situation we are considering, only one term or
the other is supposed to occur.

The point just to the right of the critical line is
of interest in  the even weight situation when $n=2m$, where the formula predicts
$$\ord_{s=m+1}L(H^{2m}(X),s)=-\dim \Hom_{ Mot_{\Z}} (\Q, H^{2m}(X)(m))$$
generalizing the pole of the Artin $L$-function ($m=0$).
As for an explicit connection to arithmetic geometry,
it is expected that
$$\Hom_{ Mot_{\Z}} (\Q, H^{2m}(X)(m)) \simeq [CH^m(X)/CH^m(X)^0]\otimes \Q$$
Of course the isomorphism should arise via a cycle map
$$CH^m(X) \ra H^{2m}(X)(m)$$
killing the cycles $CH^m(X)^0$ homologically equivalent to zero.

As we move to the left, we encounter the point $m=(n+1)/2$ for $n$ odd ($n=2m-1$),
and the conjecture predicts the order of vanishing
at the central critical point:
$$\ord_{s=m}L(H^{2m-1}(X),s)=\dim Ext^1_{Mot_{\Z}}(\Q, H^{2m-1}(X)(m))$$
It is then expected that
$$Ext^1_{Mot_{\Z}}(\Q, H^{2m-1}(X)(m))\simeq CH^m(X)^0 \otimes \Q$$
The map from cycles to extensions generalizes Kummer theory:
given a representative $Z$ for a class in $CH^m(X)^0$, we get an
exact sequence
$$0\ra H^{2m-1}(X)(m) \ra H^{2m-1}(X\setminus Z)(m) $$
$$\stackrel{\d}{\ra}
H^{2m}_Z(X)(m) \ra H^{2m}(X)(m)$$
There is a local cycle class $$cl(Z)\in
H^{2m}_Z(X)(m)$$ that maps to zero in $H^{2m}(X)(m)$,
giving rise to the desired extension:
$$0 \ra H^{2m-1}(X)(m)\ra \d^{-1}(cl(Z))\ra \Q \ra 0$$
These two classical points, central critical, $n+1-m=m=(n+1)/2$ for $n$ odd,
 and just right of it, $n+1-m=n/2+1,$ for $n$ even,
are somewhat exceptional. As $n+1-m$ moves further left ($m>n/2+1$),
one expects
$$Ext^1_{Mot_{\Z}}(\Q, H^n(X)(m))=H^{n+1}_{M,\Z}(X, \Q(m))$$
with the last group, often referred to as
{\em motivic cohomology}, defined using $K$-theory :
$$Im [ (K_{2m-n-1}(\cX))^{(m)} \ra (K_{2m-n-1}(X))^{(m)}]$$
($\cX$ is a proper flat regular $\Z$-model for $X$)
or Bloch's higher Chow groups \cite{bloch2}:
$$Im[CH^{n+1}(\cX, 2m-n-1)\otimes \Q \ra CH^{n+1}(X, 2m-n-1)\otimes \Q]$$
The latter interpretation, carrying with it the hope of representing motivic cohomology
classes quite explicitly, is more popular lately \cite{DS}.
In fact, when $m>n/2+1$, the conjectured functional equation implies
$$\ord_{s=n+1-m}L(H^n(X),s)=\dim Ext^1_{MHS^{\R}_{\R}}(\R, H^n_B(X)(m)\otimes \R)$$
where the extension occurs inside the category of
real mixed Hodge structures defined over $\R$.
So the statement on the order of vanishing follows from the
conjecture that the Hodge realization functor
induces an isomorphism
$$Ext^1_{Mot_{\Z}}(\Q, H^n(X)(m))\otimes \R \simeq Ext^1_{MHS^{\R}_{\R}}(\R, H^n_B(X)(m)\otimes \R)$$
It has been emphasized by Deligne that the regulator map (discussed below)
{\em is } this realization functor.

Our outline thus far should already make clear that the conceptual structure of the conjectures
 falls into two parts:
\bq
(1) Relation between $L$ functions and $Ext$  and $Hom$ groups in
category of motives;

(2) geometric interpretation of $Ext$  and $Hom$ groups.\eq
That is to say, in addition to the difficult problem of drawing  the lower
edge of the following triangle
$$
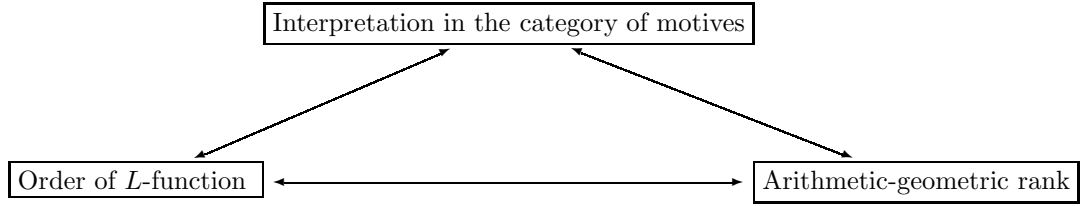
\begin{diagram}
 & &\fbox{\mbox{Interpretation in the category of motives}}& & \\
  & \ldCorresponds & & \rdCorresponds & \\
 \fbox{\mbox{Order of $L$-function} }& & \rCorresponds & & \fbox{\mbox{Arithmetic-geometric rank}}
 \end{diagram}$$
 the further problem has been created of constructing a category
 that realizes the upper vertex.

There is  a computation, convenient in practice, of
the real $Ext$ group via Deligne cohomology:
$$Ext^1_{MHS^{\R}_{\R}}(\R, H^n_B(X)(m)\otimes \R) \simeq H^{n+1}_D(X_{\R}, \R(m))$$
Using the explicit nature of Deligne
cohomology, one can construct regulator maps
$$H^{n+1}_{M,\Z}(X, \Q(m))\ra Ext^1_{MHS^{\R}_{\R}}(\R, H^n_B(X)(m)\otimes \R)$$
that can be studied independently of
a category of motives.
For example, in essentially all known cases of the Beilinson conjectures (e.g. \cite{beilinson2},\cite{deninger1}), one
constructs subgroups
$$L\subset H^{n+1}_{M,\Z}(X, \Q(m)),$$
 conjecturally  of full rank,
and studies their images in Deligne cohomology.

As we move on to the conjectures on values,
we
start with the central critical ones due to Bloch and Beilinson. Thus, $m=n+1-m=(n+1)/2$   (and $n=2m-1$).
We then have an isomorphism
$$F^mH^{2m-1}_{DR}(X)\otimes \R \simeq [ H^{2m-1}_B(X)(m-1)]^{(-1)^{m-1}}\otimes \R$$
which gives rise to the {\em period isomorphism}
$$p=p(H^{2m-1}(X,m)): \det (F^mH^{2m-1}_{DR}(X))\otimes\R \simeq \det[[ H^{2m-1}_B(X)(m-1)]^{(-1)^{m-1}}]\otimes \R$$
of real vector spaces of dimension one, where the $\det$ refers to top exterior powers. That is to say,
$$\det (F^mH^{2m-1}_{DR}(X))$$
and
$$\det[[ H^{2m-1}_B(X)(m-1)]^{(-1)^{m-1}}]$$
are viewed as two $\Q$-lines sitting inside the same real line.
An additional transcendental contribution comes from the {\em height pairing},
conjectured to be non-degenerate:
$$CH^m(X)^0\times CH^{\dim(X)+1-m}(X)^0 \ra \R$$
whose determinant gives us a regulator
$$r=r(H^{2m-1}(X)(m)) \in \R^*/\Q^*$$
Recall that conjecturally
$$d_m:=ord_{s=m}L(H^{2m-1}(X),s)=\dim CH^m(X)^0\otimes \Q$$
For any motive $M$, denote by
$$L^*(M, m)$$
the leading coefficient of $L(M,s)$ at $s=m$.
So
$$L^*(H^{2m-1}(X),m)=\lim_{s\ra m}(s-m)^{-d_m}L(H^{2m-1}(X),s)$$
The general conjecture on central critical values then says
$$r\cdot p[ \det (F^mH^{2m-1}_{DR}(X))]=L^*(H^{2m-1}(X),m)\det[[ H^{2m-1}_B(X)(m-1)]^{(-1)^{m-1}}]$$
inside $\det[[ H^{2m-1}_B(X)(m-1)]^{(-1)^{m-1}}]\otimes \R$.

Moving left to  the values at
$n+1-m \leq n/2$
($m\geq  n/2+1 $) we point out first that the period isomorphism is replaced by   an exact sequence:
$$0\ra F^mH^n_{DR}(X)\otimes \R \ra [H^n_B(X)(m-1)]^{(-1)^{m-1}}\otimes \R $$
$$\ra Ext^1_{MHS^{\R}_{\R}}(\R, H^n_B(X)(m)\otimes \R )\ra 0$$
Thus, the transcendental part should incorporate
a $\Q$-structure on
$$Ext^1_{MHS^{\R}_{\R}}(\R, H^n_B(X)(m)\otimes \R).$$
We will skip for a moment the classically interesting point $n+1-m=n/2$ for $n$ even and
assume $n+1-m<n/2$. Therefore, $m> n/2+1$ lies in the region of convergence of the $L$-function.
Then the $\Q$ structure is expected to come from the conjectured {\em regulator isomorphism}:
$$H^{n+1}_{M,\Z}(X, \Q(m)))\otimes \R \simeq Ext^1_{MHS^{\R}_{\R}}(\R, H^n_B(X)(m)\otimes \R)$$
We are led thus to an isomorphism
$$c=c(H^n(X)(m)): [\det(H^{n+1}_{M,\Z}(X, \Q(m)))]\otimes \R \simeq [\det F^mH^n_{DR}(X)]^{-1}\otimes \det
([H^n_B(X)(m-1)]^{(-1)^{m-1}})\otimes \R  $$
and
Beilinson's conjecture says
$$c[ \det(H^{n+1}_{M,\Z}(X, \Q(m)))]=L^*(H^n(X),n+1-m)[\det F^mH^n_{DR}(X)]^{-1}\otimes \det
([H^n_B(X)(m-1)]^{(-1)^{m-1}})$$
Finally, we return to the value at
$n+1-m=n/2$ ($m=n/2+1$) for $n=2m-2$ even. Here, the conjecture is identical to
the situation further left except the regulator  involves
maps  from both motivic cohomology
$$H^{n+1}_{M,\Z}(X, \Q(m))$$
and
$$CH^{m-1}(X)$$
in a manner similar to how the central critical value
incorporates periods and a height pairing. That is,
we have a map
$$[CH^{m-1}(X)/CH^{m-1}(X)^0]\ra [H^{2m-2}_B(X)(m-1)]^{(-1)^{m-1}}$$
coming from the cycle map that induces an injection
$$[CH^{m-1}(X)/CH^{m-1}(X)^0]\hra  Ext^1_{MHS^{\R}_{\R}}(\R, H^{n}_B(X)(m)\otimes \R )$$
via the quotient map
$$[H^n_B(X)(m-1)]^{(-1)^{m-1}}\otimes \R \ra Ext^1_{MHS^{\R}_{\R}}(\R, H^n_B(X)(m)\otimes \R )$$
We remark then only that  the conjectured isomorphism  is
$$H^{n+1}_{M,\Z}(X,m)\otimes \R \oplus [CH^{m-1}(X)/CH^{m-1}(X)^0]\otimes \R$$
$$\simeq Ext^1_{MHS^{\R}_{\R}}(\R, H^n_B(X)(m)\otimes \R )$$
The point is that in order to give the correct order of zero (dimension of motivic cohomology)
to the left of the critical strip, the poles to the right of the critical
strip must be canceled out.

In the Bloch-Kato conjectures, isomorphisms are
normalized more carefully,
comparing certain integral structures
one prime at a time. More precisely, the conjecture concerns the
$p$-adic valuation of the `rational part' of the
$L$-function for each prime $p$.
 We give here just a flavor of  the conjecture
by discussing a range of cases that doesn't involve too many definitions.
There is a theory of determinants of perfect complexes over
a principal ring $R$ that goes with the conjectures \cite{FP}, where the useful
facts are:
\bq
(1) If $M$ is a finitely generated free module, then $\det(M)$ is the top
exterior power.

(2) If $K$ is the field of fractions of $R$, then
for any finitely generated module $M$, the map
$M\ra M\otimes K$ induces
$$\det(M)\otimes_R K\simeq \det(M\otimes_RK)$$
canonically.

(3) If $M$ is finitely generated and $M'$ is its torsion-free
quotient, then the natural isomorphism
$$\det(M)\otimes K\simeq \det(M')\otimes K$$
takes  a generator of $\det(M)$ to $1/r$ times a generator of $\det(M')$,
where $r$ is a characteristic element for the torsion submodule of $M$.
\eq
We note in regard to these facts that the determinant module is
always free of rank one, so it is only the morphisms
that are the relevant data.
Now assume that $m>\mbox{min}\{n, \dim(X)\}$ so that $F^mH^n_{DR}=0$ and
$$[H^n_B(X)(m-1)^{(-1)^{m-1}}\otimes \R] \simeq Ext^1_{MHS^{\R}_{\R}}(\R, H^n_B(X)(m)\otimes \R )$$
Therefore, we should have
$$H^{n+1}_{M,\Z}(X, \Q(m))\otimes \R
\simeq Ext^1_{MHS^{\R}_{\R}}(\R, H^n_B(X)(m)\otimes \R)$$
$$\simeq [H^n(X(\C), \Z(m-1))^{(-1)^{m-1}}\otimes \R]$$
Let $\g$ be a generator of $$\det(H^n(X(\C), \Z(m-1))^{(-1)^{m-1}})$$
and let $$\om_{\g}\in \det (H^{n+1}_{M,\Z}(X, \Q(m)))$$
be an element mapping to $$L^*(H^n(X)(m),n+1-m)\g.$$ As before,
let $S$ be a set of primes including those of bad reduction for $X$, the Archimedean prime,
and a fixed prime $p$.  We denote by $G_S$ the Galois group of the
maximal extension of $\Q$ unramified outside the primes in $S$.
There is a Chern-class map
$$H^{n+1}_{M,\Z}(X, \Q(m))\ra H^1(G_S, H^n_p(X)(m))$$
that  conjecturally induces an isomorphism
$$H^{n+1}_{M,\Z}(X, \Q(m))\otimes \Q_p \simeq H^1(G_S, H^n_p(X)(m))$$
That is to say, there are {\em two} maps
$$\begin{diagram}
H^{n+1}_{M,\Z}(X, \Q(m))
&\rTo & Ext^1_{MHS^{\R}_{\R}}(\R, H^n_B(X)(m)\otimes \R)\\
 \dTo& & \\
 H^1(G_S, H^n_p(X)(m)) & &
 \end{diagram}$$
 that are both supposed to induce isomorphisms upon changing coefficients,
 and with which the rational part will be controlled.
Let $$z_{\g}\in \det (H^1(G_S, H^n_p(X)(m)))$$ be the image of $\om_{\g}$
so that  we have the diagram
$$\begin{diagram}
\om_{\g} & \rMapsto& L^*(H^n(X)(m), n+1-m)\g \\
 \dMapsto& & \\
 z_g & & \end{diagram}
 $$
Then the conjecture is that
$$[\det(H^1(G_S, H^1 (G, H^n_p(X)(m))): \Z_p z_{\g}]=|H^0(G_S,H^1_f(G, H^n_p(X)(m)))||H^2(G_S,H^1_f(G, H^n_p(X)(m)))|$$
It is easy to see that the $H^0$ term is finite, but the finiteness of
the $H^2$, like that of $Sha$ in the case of elliptic curves,  must also be conjectured.
In the uniform formalism, the desired
equality  is interpreted as the assertion that  $z_{\g}$ is a $\Z_p$-basis of
the determinant of the perfect $\Z_p$-complex
$R\G(G_S, H^n_p(X)(m))$.

Extraction of the rational part is supposed to lead eventually
to a
{\em $p$-adic} $L$-function
$${\cal L}^{(p)}(H^n(X))$$ that exercises control over Galois cohomology (i.e.,
Selmer groups) and
Diophantine invariants.
This $p$-adic theory appears so far to be the best strategy  for applying  the theory of $L$-functions to the
elucidation of Diophantine structures (\cite{rubin}, \cite{kato2}, \cite{SU}).
\section{Remark}

We conclude with the warning that there is a conspicuous deficiency in theory of motives:
This is that even in the best of possible worlds (ours), only {\em abelian} invariants
are accessible, such as
$$CH^m(X).$$
These abelian invariants do not yield in general information about
$$X(\Q)$$
 and leave thereby untouched the most
basic questions of Diophantine geometry.
This is an artifact of the fact that the theory of motives as presently
developed  is implicitly modeled  on the
theory of abelian varieties and $H_1$. Attempts to redress this deficiency for certain varieties are contained in
$$\mbox{Grothendieck's {\em anabelian program}}$$
(\cite{grothendieck}, \cite{NTM})
that concerns itself with the theory of {\em pro-finite $\pi_1$'s}. The technology of motives ends up
contributing here as well because the Diophantine aspect of this theory \cite{kim}
assigns an interesting role
to  {\em motivic fundamental groups} \cite{deligne2}, where
$Ext$ groups are replaced by
$$\mbox{\em classifying spaces
for non-abelian torsors} $$ However, what is entirely missing as yet is
an analogue of the $L$-function.

\end{document}